\documentclass[tikz,12pt]{article}
\usepackage{amsmath,amssymb,amsthm,graphicx}
\usepackage{geometry}
\usepackage{algorithm}
\usepackage{algpseudocode}
\usepackage{float}
\usepackage{caption}

\usepackage{tikz}

\newcommand{\pbb}{{\mathbf p}}
\newcommand{\vbb}{{\mathbf v}}
\newcommand{\R}{\mathbb{R}}
\newtheorem{remark}{Remark}

\DeclareMathOperator{\tr}{tr}

\algrenewcommand\algorithmicindent{1.2em}

\title{Optimal placement and tuning of pointwise dampers
for vibrating strings via a Lyapunov framework}

\author{
Josip Tamba\v{c}a\\
Department of Mathematics\\
Faculty of Science, University of Zagreb\\
Bijeni\v{c}ka cesta 30, 10000 Zagreb, Croatia\\
\texttt{tambaca@math.hr}
\and
Ninoslav Truhar \\ 
School of Applied Mathematics and Informatics\\
Josip Juraj Strossmayer University of Osijek\\
Trg Ljudevita Gaja 6, 31000 Osijek, Croatia\\
\texttt{ntruhar@mathos.hr}
}

\date{}

\begin{document}

\maketitle

\begin{abstract}
We study the optimal placement and tuning of a small number of pointwise viscous dampers for a vibrating string. Starting from a finite element discretization of the damped wave equation, the system is transformed into a first-order phase-space formulation, which enables a unified Lyapunov trace framework.

Three optimization criteria are considered: average total energy, average total displacement, and energy for a fixed initial state. For all criteria, explicit gradient formulas with respect to damper positions and viscosities are derived, requiring only one primal and one dual Lyapunov solve.

Due to the strong non-convexity of the problem, a simple heuristic based on an explicit single-damper formula is proposed to generate effective initial guesses. Numerical examples illustrate the influence of spectral selection and discretization on the optimal damping configuration.
\end{abstract}



\section{Introduction and Preliminaries}
\label{sec:introduction_preliminaries}

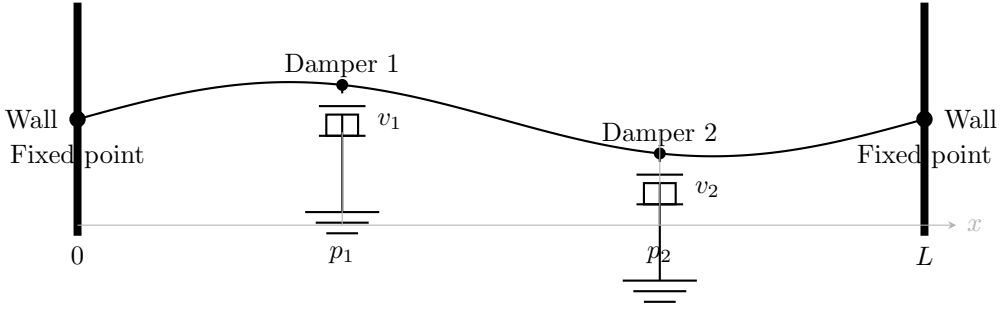
\begin{figure}[h!]
\centering
\begin{tikzpicture}[scale=1.4, every node/.style={font=\footnotesize}]
  \def\L{8}            
  \def\A{0.35}         
  \def\xone{2.5}       
  \def\xtwo{5.5}       

  \pgfmathsetmacro{\yone}{\A*sin(360*\xone/\L)}
  \pgfmathsetmacro{\ytwo}{\A*sin(360*\xtwo/\L)}

  \tikzset{
    wall/.style   ={line width=3pt},
    string/.style ={thick},
    damper/.style ={thick},
  }

  \draw[string, smooth, domain=0:\L, samples=200]
    plot (\x, {\A*sin(360*\x/\L)});

  \draw[wall] (0,-1.1) -- (0,1.1) node[midway, left=2pt] {Wall};
  \draw[wall] (\L,-1.1) -- (\L,1.1) node[midway, right=2pt] {Wall};

  \fill (0,0) circle (2.2pt) node[below=4pt] {\strut Fixed point};
  \fill (\L,0) circle (2.2pt) node[below=4pt] {\strut Fixed point};

  \fill (\xone,\yone) circle (1.6pt);
  \fill (\xtwo,\ytwo) circle (1.6pt);

  \draw[damper] (\xone,\yone) -- (\xone,\yone-0.08);
  \draw[damper] (\xone-0.22,\yone-0.20) -- (\xone+0.22,\yone-0.20);
  \draw[damper] (\xone-0.15,\yone-0.48) rectangle (\xone+0.15,\yone-0.28);
  \draw[damper] (\xone,\yone-0.28) -- (\xone,\yone-0.48);
  \draw[damper] (\xone-0.22,\yone-0.48) -- (\xone+0.22,\yone-0.48);
  \draw[damper] (\xone,\yone-0.48) -- (\xone,\yone-1.20);
  \draw[damper] (\xone-0.35,\yone-1.20) -- (\xone+0.35,\yone-1.20);
  \draw[damper] (\xone-0.25,\yone-1.30) -- (\xone+0.25,\yone-1.30);
  \draw[damper] (\xone-0.15,\yone-1.40) -- (\xone+0.15,\yone-1.40);
  \node at (\xone,\yone+0.18) {\strut Damper 1};
  \node at (\xone+0.45,\yone-0.34) {$v_1$};

  \draw[damper] (\xtwo,\ytwo) -- (\xtwo,\ytwo-0.08);
  \draw[damper] (\xtwo-0.22,\ytwo-0.20) -- (\xtwo+0.22,\ytwo-0.20);
  \draw[damper] (\xtwo-0.15,\ytwo-0.48) rectangle (\xtwo+0.15,\ytwo-0.28);
  \draw[damper] (\xtwo,\ytwo-0.28) -- (\xtwo,\ytwo-0.48);
  \draw[damper] (\xtwo-0.22,\ytwo-0.48) -- (\xtwo+0.22,\ytwo-0.48);
  \draw[damper] (\xtwo,\ytwo-0.48) -- (\xtwo,\ytwo-1.20);
  \draw[damper] (\xtwo-0.35,\ytwo-1.20) -- (\xtwo+0.35,\ytwo-1.20);
  \draw[damper] (\xtwo-0.25,\ytwo-1.30) -- (\xtwo+0.25,\ytwo-1.30);
  \draw[damper] (\xtwo-0.15,\ytwo-1.40) -- (\xtwo+0.15,\ytwo-1.40);
  \node at (\xtwo,\ytwo+0.18) {\strut Damper 2};
  \node at (\xtwo+0.45,\ytwo-0.34) {$v_2$};

%
  \draw[gray!60,-stealth] (0,-1.0) -- (\L+0.3,-1.0) node[right] {$x$};

  \draw[gray!60] (\xone,-1.0) -- (\xone,-0.12);
  \draw[gray!60] (\xtwo,-1.0) -- (\xtwo,-0.12);

  \node[below=4pt] at (\xone,-1.0) {$p_1$};
  \node[below=4pt] at (\xtwo,-1.0) {$p_2$};

  \node[below=4pt] at (0,-1.0) {$0$};
  \node[below=4pt] at (\L,-1.0) {$L$};

\end{tikzpicture}
\caption{Vibrating string with two grounded (pointwise) dampers at positions $p_1$ and $p_2$, with viscosities $v_1$ and $v_2$.}
\label{fig:string_dampers_grounded}
\end{figure}

The primary objective of this paper is to determine the optimal positions and viscosities of a small number $r$ of pointwise dampers, with $r \ll n$, for the free oscillations of a string. We consider a vibrating string fixed at both ends and equipped with pointwise viscous dampers, as illustrated in Figure~\ref{fig:string_dampers_grounded}. This configuration represents the central problem studied throughout the paper, rather than a mere motivating example.

Recent work on deformable strings \cite{Giordano2024} and on wave and beam
equations with localized Kelvin--Voigt damping \cite{LiuZhang2025,Zhang2024}
highlights the importance of localized damping mechanisms in one-dimensional
systems. In finite-dimensional mechanical systems, the selection of damping parameters
has been extensively studied; in particular, the recent work \cite{LiTisseur2026} shares
several key ideas with our approach, as it adopts a similar theoretical framework and considers
the optimization of damping coefficients (viscosities) for a fixed damper configuration.
Moreover, a substantial part of the literature discussed in \cite{LiTisseur2026} overlaps with
the references relevant to the present work, and therefore we do not repeat those citations here.

For that purpose, we consider the damped wave equation
\begin{equation}\label{strong}
    \rho A \frac{\partial^2 u(x,t)}{\partial t^2}
    + c \frac{\partial u(x,t)}{\partial t}
    - \frac{\partial }{\partial x}\left( T\frac{\partial u(x,t)}{\partial x}\right) = 0,
\end{equation}
on the space domain $[0,L]$, where $u(x,t)$ denotes the transverse displacement of the string, $T$ is the tension, $\rho$ is the volume mass density, $A$ is the area of the cross--section, and $c$ is the damping coefficient. The string is clamped at both ends, which leads to the boundary conditions
\begin{equation}\label{bc}
    u(0,t) = 0, \qquad u(L,t) = 0.
\end{equation}
Note that in this setting we allow that $\rho$ and $T$ are $L^\infty(0,L)$ functions which are uniformly positive, i.e., there are numbers $\alpha_{\rho A}>0$ and $\alpha_{T}>0$ such that
\begin{equation}\label{positivity}
\inf{\rho A ([0,L])} > \alpha_{\rho A}, \qquad \inf{T([0,L])} > \alpha_{T}.
\end{equation}
This includes problems with nonconstant cross--section of the string or non-constant tension.
To model the pointwise dampers this setting also includes the damping coefficient which is a sum of pointwise Dirac measures.
For instance localized
damping is introduced through dampers located at positions \( p_1 \) and \( p_2 \), where the damping forces are given by:
\begin{equation*}
c_1 \frac{\partial u(p_1, t)}{\partial t} = \text{Force at damper 1}, \quad c_2 \frac{\partial u(p_2, t)}{\partial t} = \text{Force at damper 2},
\end{equation*}
that corresponds to $c=c_1 \delta_{p_1} + c_2 \delta_{p_2}$, where $\delta_{p_i}$ is the Dirac measure positioned at point $p_i$. In general we assume that there are $r$ dampers at positions given by points $p_k \in [0,L]$ and viscosities $v_k$, for $k=1,\ldots,r$ and that there is internal constant damping given by a function $c_0\in L^\infty(0,L)$. Then
$$
c=c_0 + \sum_{k=1}^r v_k \delta_{p_k}.
$$

To derive the weak form of the governing equation, we multiply the partial differential equation by a smooth enough test function \( v(x) \) and integrate over the domain \( [0, L] \):
\begin{equation}\label{weak}
\int_0^L \rho A \frac{\partial^2 u}{\partial t^2} v \, dx +  \int_0^L c  \frac{\partial u}{\partial t} v\, dx
-  \int_0^L \frac{\partial }{\partial x}\left( T\frac{\partial u(x,t)}{\partial x}\right) v \, dx = 0.
\end{equation}


\subsection*{Finite Element Formulation}

To approximate the solution of the problem \eqref{weak} numerically, we follow the standard FEM practice~\cite{Zienkiewicz2013}.
Let us first split the space domain $[0,L]$ into $N$ subintervals by points
$$
0= x_0 < x_1 < \cdots < x_{N-1} < x_N = L.
$$
We use linear Lagrange elements, which are simple 1D elements where the displacement $u(x)$ is approximated linearly within each element. For an element spanning from $x_i$ to $x_{i+1}$, the shape functions $N_1(x)$ and $N_2(x)$ are defined as:
\begin{equation*}
N_1(x) = \frac{x_{i+1} - x}{x_{i+1} - x_i}, \quad N_2(x) = \frac{x - x_i}{x_{i+1} - x_i}.
\end{equation*}
Based on these for $i\in \{1,\ldots, N-1\}$ we define functions
$$
\varphi_i (x) = \left\{\begin{array}{ll}
\frac{x-x_{i-1}}{x_{i} - x_{i-1}}, &  x \in [x_{i-1},x_i],\\
\frac{x_{i+1} - x}{x_{i+1} - x_{i}}, &  x \in [x_{i},x_{i+1}],\\
0,& \mbox{otherwise}
\end{array}
\right.
$$
that form the basis of FEM, $L\{\varphi_1, \ldots, \varphi_{N-1}\}$. Then we look for the solution of the problem \eqref{weak} in the form
$$
u(x, t) \approx \sum_{j=1}^{N-1} u_j(t) \varphi_j(x),
$$
where $u_j(t)$ are the nodal displacements at the element nodes.

\subsection*{System Matrices}
Substituting the finite element approximation into the weak form, we derive the global system matrices.

\paragraph{Mass Matrix ($M$):}
The elements of the global mass matrix are defined as
\begin{equation}
M_{ij} = \int_0^L \rho A\, \varphi_i(x) \varphi_j(x)\, dx.
\end{equation}
For a linear element of length $h$ with shape functions in barycentric coordinates
$N_1(\xi)=1-\xi$, $N_2(\xi)=\xi$, and for constant cross--section and density, the local (elementary) mass matrix becomes
\begin{equation}
M^{(e)} = \rho A h \int_0^1
\begin{bmatrix}
(1-\xi)^2 & (1-\xi)\xi \\[6pt]
\xi(1-\xi) & \xi^2
\end{bmatrix} d\xi
= \frac{\rho A h}{6}
\begin{bmatrix}
2 & 1 \\[6pt]
1 & 2
\end{bmatrix}.
\end{equation}
In component form, the entries are
\begin{equation}
M_{ij}^{(e)} = \frac{\rho A h}{6}
\begin{cases}
2, & i=j, \\[6pt]
1, & i \neq j,
\end{cases}
\qquad i,j=1,2.
\end{equation}

Assembling all element contributions on a uniform mesh with $N$ elements
of length $h=L/N$ leads to the global (consistent) mass matrix
\begin{equation}
\mathbf M = \frac{\rho A h}{6}
\begin{bmatrix}
4 & 1 &   &   &   \\
1 & 4 & 1 &   &   \\
  & 1 & 4 & 1 &   \\
  &   & \ddots & \ddots & \ddots \\
  &   &        & 1 & 4 & 1 \\
  &   &        &   & 1 & 4
\end{bmatrix}_{(N-1)\times(N-1)}.
\label{eq:global-M}
\end{equation}
Here, the Dirichlet boundary conditions $u(0)=u(L)=0$ are already included in the function space. Therefore for a constant function $\rho A$ we have
\begin{equation} \label{eq:reduced-M}
\mathbf M_{\text{red}} = \frac{\rho A h}{6}\,
\operatorname{tridiag}(1,\,4,\,1)
= \frac{\rho A h}{3}\,\operatorname{tridiag}\!\Big(\tfrac12,\,2,\,\tfrac12\Big),
\end{equation}
which is symmetric, tridiagonal, and positive definite.


\paragraph{Stiffness Matrix ($K$):}
The elements of the global stiffness matrix are computed as:
\begin{equation}
K_{ij} = \int_0^L T \frac{\partial \varphi_i}{\partial x} \frac{\partial \varphi_j}{\partial x} \, dx, \quad i,j\in \{1,\ldots,N-1\}.
\end{equation}
For a linear element of length $h$, and constant tension, the local stiffness matrix is:
\begin{equation}
K^{(e)} = \frac{T}{h}
\begin{bmatrix}
1 & -1 \\[6pt]
-1 & 1
\end{bmatrix}.
\end{equation}
Assembling over $N$ uniform elements of length $h=L/N$ gives the global stiffness matrix
\begin{equation}
\mathbf K = \frac{T}{h}
\begin{bmatrix}
2 & -1 &   &   &   \\
-1 & 2 & -1 &   &   \\
   & -1 & 2 & -1 &   \\
   &    & \ddots & \ddots & \ddots \\
   &    &        & -1 & 2 & -1 \\
   &    &        &    & -1 & 2
\end{bmatrix}_{(N-1)\times(N-1)}.
\label{eq:global-K}
\end{equation}
Again, the Dirichlet boundary conditions $u(0)=u(L)=0$ are already included in the function space of FEM. Thus, for a constant tension we have
\begin{equation} \label{eq:reduced-K}
\mathbf K_{\text{red}} = \frac{T}{h}\,
\operatorname{tridiag}(-1,\,2,\,-1).
\end{equation}

\paragraph{Damping Matrix ($D$):}

The total damping matrix consists of localized pointwise dampers and internal
(material) damping. Its elements are defined by
\begin{equation}\label{dij}
D_{ij} = \int_0^L c \varphi_i(x) \varphi_j(x) dx = \int_0^L c_0 \varphi_i(x) \varphi_j(x) dx + \sum_{k=1}^r v_k \varphi_i(p_k) \varphi_j(p_k)
\end{equation}
In the sequel the matrix for the internal damping we denote by $D_0$. Pointwise dampers we resolve in the following way.

Let us define vector functions $\hat{d}:[0,L] \to \R^{N-1}$ by
\begin{equation}
\hat{d} (p) = \left[\begin{array}{c} \varphi_1(p)\\ \varphi_2(p)\\ \vdots \\ \varphi_{N-1}(p)
\end{array}\right], \qquad p\in[0,L].
\label{eq:def_hatd}
\end{equation}
Then
\begin{equation}
D = D_{0} + \sum_{k=1}^r v_k {\hat{d}(p_k) \hat{d}(p_k)^T}.
\label{eq:D-total}
\end{equation}

{\begin{remark}\label{rsmooth}\em
By definition the component functions of $\hat{d}$ are piecewise linear functions on each subinterval determined by the FEM mesh and globally continuous. Thus these functions are smooth up to the finite points of the FEM mesh, $\{x_1, \ldots, x_{N-1}\}$. Their derivatives can be easily calculated
$$
\hat{d}' (p) = \left[\begin{array}{c} \varphi_1'(p)\\ \varphi_2'(p)\\ \vdots \\ \varphi_{N-1}'(p)
\end{array}\right],
$$
where
$$
\hat{d}_j' (p) = \varphi_j'(p) = \left\{\begin{array}{ll}
\frac{1}{h} & p \in [x_{j-1},x_j]\\
-\frac{1}{h} & p \in [x_{j},x_{j+1}]\\
0 & \mbox{otherwise}
\end{array}
\right..
$$
\end{remark}

It remains to describe the internal (material) damping, which we model either
by the classical Rayleigh form
\[
D_{0} = \alpha_D M + \beta_D K ,
\]
or by assigning a small fraction $\zeta$ of critical damping, using the modal
expression
\[
D_{0} = \zeta \,
M^{-1/2}\sqrt{\,M^{1/2} K M^{1/2}\,}\, M^{-1/2},
\]
as used in structural vibration modelling
(see \cite{Zienkiewicz2013, KuzmTomljTruh12}).

\medskip

Collecting all contributions described above, the FEM model leads to the
second-order matrix equation
\begin{equation}
M \ddot{x}(t) + D \dot{x}(t) + K x(t) = 0,
\qquad
x(0)=x_0,\quad \dot{x}(0)=\dot{x}_0.
\label{eq:SecondOrderFEM}
\end{equation}
All matrices $M,D,$ and $K$ are symmetric, and due to uniform positivity \eqref{positivity} they are positive definite.

It should be emphasized that matrices $M,K$ are in \eqref{eq:global-M}, \eqref{eq:reduced-M}, \eqref{eq:global-K} and \eqref{eq:reduced-K}
for constant coefficients. However, in numerical examples we consider models with
non-constant coefficients. In that case matrices $M$ and $K$ are obtained using numerical
integration formula (trapezoidal rule).

\section{Model Transformation and Optimization Framework}

In what follows, we set $n = N-1$. For the purpose of damping optimization, we follow an approach similar to that developed in several earlier papers on damping optimization (positions and viscosity) \cite{TRUHVES09, KuzmTomljTruh12, TRUTOMVES2014}.
More recently, adaptive and model-reduction-based optimization strategies have been proposed in \cite{PrzybillaUgricaTruharBenner2025}, where the damping positions are optimized through controllability-based decoupling techniques.

\medskip

Following the classical methodology for linear vibrational systems
\cite{MullerSchiehlen85, VES2011}, we now transform the finite element
system \eqref{eq:SecondOrderFEM} into an equivalent first--order
representation.
Let $\Phi$ be a non-singular matrix that simultaneously diagonalizes the pair
$(M,K)$:
\begin{equation*}
\Phi^T M \Phi = I,
\qquad
\Phi^T K \Phi = \Omega^{2},
\end{equation*}
where $\Omega = \mathrm{diag}(\omega_1,\ldots,\omega_n)$ contains the natural
frequencies of the undamped system.

Define the phase-space variables
\begin{equation}\label{Def:PhaseSpacevec}
\mathbf{y}_1(t) = \Omega \Phi^{-1} \mathbf{x}(t),
\qquad
\mathbf{y}_2(t) = \Phi^{-1} \dot{\mathbf{x}}(t),
\end{equation}
and set
\[
\mathbf{y}(t) =
\begin{bmatrix}
\mathbf{y}_1(t) \\[4pt]
\mathbf{y}_2(t)
\end{bmatrix}.
\]

Then the second-order system
\[
M \ddot{x}(t) + D \dot{x}(t) + K x(t) = 0
\]
is equivalent to the first-order system
\begin{equation}\label{eq:PhaseSpaceSystem}
\dot{\mathbf{y}}(t) = A \mathbf{y}(t),
\qquad
A \doteq A(\mathbf{p},\mathbf{v}) =
\begin{bmatrix}
0 & \Omega \\
-\Omega & -C_{\Phi}
\end{bmatrix},
\end{equation}
where
\[
C_{\Phi} \doteq  C_{\Phi}(\mathbf{p},\mathbf{v}) = \Phi^T D \Phi = \sum_{k=1}^r v_k {\Phi^T} {\hat{d}(p_k) \hat{d}(p_k)}^T {\Phi}\,
\]
where $\mathbf{p}=(p_1,\dots,p_r)$ and $\mathbf{v}=(v_1,\dots,v_r)$
denote the vectors of damper positions and corresponding viscosities,
respectively, and $\hat{d}(p_k)$ is the $n$-dimensional vector obtained
from the FEM formulation described in ``Damping Matrix (D)'' as in
\eqref{eq:def_hatd}.


This reformulation places the dynamics in the framework suitable for
Lyapunov-based damping optimization, following the theory developed in
\cite{Nakic02, Brabender98, TRUHVES09, TRUTOMVES2014}.
In particular, it allows us to define optimization criteria through the trace of
solutions to Lyapunov equations, which quantify the mean energy or mean
displacement over all initial conditions.

\section{Optimization Criteria}

Following the phase--space formulation introduced in \eqref{Def:PhaseSpacevec}, we now describe the
criteria used to optimize the damping distribution.
We consider three measures of system performance: the average total energy, the average total displacements, and the energy corresponding to
fixed initial conditions.
All three lead to trace minimization problems involving solutions of Lyapunov equations.

\subsection{Average Total Energy}
\label{subsec_average_total_energy}

The Euclidean norm of the phase--space vector $\mathbf{y}(t)$ is directly related
to the total energy of the vibrating system.
From the definition of the phase variables in \eqref{Def:PhaseSpacevec}, we obtain
\begin{equation*}
\mathbf{y}(t)^T \mathbf{y}(t)
= \|\mathbf{y}_1(t)\|^2 + \|\mathbf{y}_2(t)\|^2
= x(t)^T K\, x(t) + \dot{x}(t)^T M\, \dot{x}(t)
= 2 E(t),
\end{equation*}
where $E(t)$ denotes the total mechanical energy at time~$t$.
This identity implies that all admissible phase--space representations are
unitarily equivalent with respect to energy, and therefore any convenient
coordinate system may be used for optimization.

For the first optimization criterion, we minimize the mean value of the total
energy over all possible initial conditions.
As shown in \cite{Nakic02, Brabender98}, this problem is equivalent to minimizing
\begin{equation}\label{Trace_Energy}
\operatorname{tr}(Z X) \;\longrightarrow\; \min,
\end{equation}
where $X$ is the unique solution of the Lyapunov equation
\begin{equation}\label{Lyapunov_Energy}
A^T X + X A = -I.
\end{equation}
The weighting matrix $Z = Z_{s_0} \oplus Z_{s_0}$ selects the eigenfrequencies to be
controlled.
If the goal is to damp the first ${s_0}$ modes
$0 < \omega_1 < \omega_2 < \dots < \omega_{s_0}$,
then
\begin{equation}\label{Z_energy}
Z = Z_{s_0} \oplus Z_{s_0},
\qquad
Z_{s_0} =
\begin{bmatrix}
I_{s_0} & \\
& 0_{(n-{s_0})}
\end{bmatrix}.
\end{equation}
Further details on constructing $Z$ may be found in \cite{Nakic02}.

An equivalent formulation is obtained by solving the dual Lyapunov equation
\begin{equation}\label{Lyapunov_Dual_Energy}
A Y + Y A^T = -Z,
\end{equation}
because
\begin{equation}\label{Trace_Dual_Energy}
\operatorname{tr}(Y) = \operatorname{tr}(Z X).
\end{equation}
Thus minimizing \eqref{Trace_Energy} is equivalent to minimizing the trace of the
solution $Y$ of \eqref{Lyapunov_Dual_Energy}.

\subsection{Average Total Displacements}
\label{subsec_average_total_displacements}

To the best of our knowledge, apart from the recent work by Truhar and Veseli\'c~\cite{TruharVeselic2025},
the criterion of minimizing the \emph{average total displacement} has not been considered in the literature.

As one can find in~\cite{TruharVeselic2025}, the mean--square displacement criterion reduces to minimizing the trace of the
solution $\widehat{X}$ of the Lyapunov equation
\begin{equation*}
A^{T}\widehat{X} + \widehat{X} A = -Z_{ds},
\qquad
Z_{ds} = \begin{bmatrix} \widehat{K}^{-1} & 0 \\[2pt] 0 & 0 \end{bmatrix},
\end{equation*}
where the system matrix $A$ is defined in \eqref{eq:PhaseSpaceSystem} and
$\widehat{K}^{-1}$ is given in
\begin{equation}\label{Def_K_Hat}
\Omega^{-1} \Phi^T \Phi \Omega^{-1} = V \Lambda_K^{-1} V^T \doteq \widehat{K}^{-1},
\end{equation}

Note  that in both criteria one needs to minimize a trace of the corresponding {Lyapunov} equation,
with the same system matrix but a different right-hand side.

In the next section, we introduce a further criterion motivated by the case when the initial conditions are exactly known,
which can again be formulated as the minimization of the trace of the solution to the same Lyapunov equation whose
system matrix is given in~\eqref{eq:PhaseSpaceSystem}.

\subsection{Energy-based optimization criterion for fixed initial conditions ($y_0^{*}Xy_0$)}
\label{subsec_Energy-based optimization}

As discussed earlier, one of the criteria for damping optimization
is the minimization of the norm of the solution of the first--order
linear system
\begin{equation}\label{eq:ode_first_order}
y'(t)=A(\mathbf p,\mathbf v)\,y(t), \qquad y(0)=y_0 ,
\end{equation}
where $y_0 \in \mathbb{R}^{2n}$ (or $\mathbb{C}^{2n}$) is a fixed initial condition.

The solution of \eqref{eq:ode_first_order} is given by
\[
y(t)=e^{A(\mathbf p,\mathbf v)t}\,y_0 .
\]
The squared Euclidean norm of the state at time $t$ satisfies
\[
\|y(t)\|^2
= y(t)^{*}y(t)
= y_0^{*} e^{A(\mathbf p,\mathbf v)^{*} t}\,
e^{A(\mathbf p,\mathbf v) t}\,y_0 .
\]

The total energy associated with the initial condition $y_0$
is defined as the time integral of the state norm,
\begin{equation}\label{eq:energy_y0}
E(\mathbf p,\mathbf v;y_0)
:=\int_{0}^{\infty}\|y(t)\|^2\,dt
=\int_{0}^{\infty}
y_0^{*} e^{A(\mathbf p,\mathbf v)^{*} t}\,
e^{A(\mathbf p,\mathbf v) t}\,y_0 \,dt .
\end{equation}

Introducing the matrix
\[
X(\mathbf p,\mathbf v)
:=\int_{0}^{\infty}
e^{A(\mathbf p,\mathbf v)^{*} t}\,
e^{A(\mathbf p,\mathbf v) t}\,dt ,
\]
the energy can be written in the compact form
\[
E(\mathbf p,\mathbf v;y_0)= y_0^{*} X(\mathbf p,\mathbf v)\, y_0 \,,
\]
where the matrix $X(\mathbf p,\mathbf v)$ satisfies the Lyapunov equation
\eqref{Lyapunov_Energy}, that is
\begin{equation*}
A(\mathbf p,\mathbf v)^{T} X + X A(\mathbf p,\mathbf v) = -I ,
\end{equation*}
where, in the real case, $A^{*}=A^{T}$.

%

Moreover,
\[
E(\mathbf p,\mathbf v;y_0)
=
\operatorname{tr}\!\bigl(X(\mathbf p,\mathbf v)\, y_0y_0^{*}\bigr).
\]

Note that this representation is of the same form as the trace objective
\[
f(\mathbf p,\mathbf v)
=
\operatorname{tr}\!\bigl(Z X(\mathbf p,\mathbf v)\bigr)
\]
from the previous section, with the particular choice
\[
Z = y_0 y_0^{*}.
\]
Accordingly, the gradient formulas derived in
Section~\ref{Sec:Gradoftracefunc} apply directly in this setting.

Following the general framework, we introduce the matrix
$Y_E = Y_E(\mathbf p,\mathbf v)$ as the unique solution of the dual Lyapunov equation
\[
A(\mathbf p,\mathbf v) Y_E + Y_E A(\mathbf p,\mathbf v)^T
= -\,y_0 y_0^{*}.
\]

Up to the authors' knowledge, the criterion of minimizing the total energy associated with a specific initial condition ($y_0$), i.e. the minimization of quantity \eqref{eq:energy_y0}, was first considered in the context of optimal damping of free vibrations in multi-degree-of-freedom systems with proportional damping \cite{LelasNakic2024}, where it was shown that the resulting optimal damping strongly depends on the shares of potential and kinetic energy in the initial energy. In \cite{LelasNakic2024}, due to the simplicity of the system under consideration, it was possible to perform a modal analysis and analytically calculate \eqref{eq:energy_y0} for any initial condition.

Building on the three Lyapunov--based criteria introduced above, the optimization tasks can now be expressed within a unified framework:
the quantities to be minimized are the traces of the solutions of Lyapunov equations whose system matrix is given in~\eqref{eq:PhaseSpaceSystem},
with the damper positions and viscosities entering only through the modified damping term.
This observation allows us to formulate a single optimization problem in which both parameter sets are treated consistently.
In the next section, we develop an iterative gradient-based optimization scheme that simultaneously updates the viscosities \(\mathbf{v}\)
and damper positions \(\mathbf{p}\) so as to minimize any of the three criteria introduced above.

\section{Gradient of the trace objective}
\label{Sec:Gradoftracefunc}

The objective function in all three optimization problems, is locally non-differentiable at finite discrete attachment points along the FEM string model
due to the choice of the $P1$ FEM basis functions. However, it is well known that
gradient-based methods exhibit remarkable empirical robustness in such environments, and successfully traverses these non-smooth points, see e.g. \cite{LO} or \cite{BS}

In gradient-based optimization methods, the computational bottleneck often
arises from evaluating derivatives of the objective function.
In our case, the objective function is
\[
f(\mathbf{p},\mathbf{v})
:= \operatorname{tr}\!\bigl(Z\,X(\mathbf{p},\mathbf{v})\bigr),
\]
where \(Z\) is a fixed symmetric matrix associated with the chosen minimization criterion
and \(X(\mathbf{p},\mathbf{v})\) is defined implicitly as the unique solution of the
Lyapunov equation
\begin{equation}
A(\mathbf{p},\mathbf{v})^{T} X + X\,A(\mathbf{p},\mathbf{v}) = -R_{hs}\, ,
\label{eq:Lyap_grad_main0}
\end{equation}
where the right-hand side matrix \(R_{hs}\) is determined by the same criterion.
Hence, the dependence of \(f\) on the parameters
\(\mathbf{p}=(p_1,\dots,p_r)\) and \(\mathbf{v}=(v_1,\dots,v_r)\)
is entirely mediated through the system matrix \(A(\mathbf{p},\mathbf{v})\).

\medskip

\paragraph{Integral representation.}
Since the matrix \(A(\mathbf{p},\mathbf{v})\) is assumed to be Hurwitz,
the solution of \eqref{eq:Lyap_grad_main0} admits the integral representation
\[
X(\mathbf{p},\mathbf{v})
= \int_0^\infty
e^{A(\mathbf{p},\mathbf{v})^{T} t}
e^{A(\mathbf{p},\mathbf{v}) t}\,\mathrm{d}t .
\]
This classical representation follows from the general theory of Lyapunov
equations; see, for instance,  \cite{Golub1996,LANROD95}.
It will play a central role in the derivation of explicit gradient formulas.

\medskip

\paragraph{Differentiation with respect to a parameter.}
Let \(\vartheta\) denote one of the parameters \(v_i\) or \(p_i\).
By linearity of the trace, differentiation of the objective function yields
\[
\frac{\partial f}{\partial \vartheta}
= \operatorname{tr}\!\left(
Z\,\frac{\partial X}{\partial \vartheta}
\right),
\]
where the derivative \(\partial X / \partial \vartheta\) is induced implicitly
through the dependence of the system matrix
\(A(\mathbf{p},\mathbf{v})\) on the parameter \(\vartheta\).

Differentiating \eqref{eq:Lyap_grad_main0} with respect to \(\vartheta\) gives
\[
A^{T} \frac{\partial X}{\partial \vartheta}
+ \frac{\partial X}{\partial \vartheta} A
 =
- \left(
\frac{\partial A^{T}}{\partial \vartheta} X
+ X \frac{\partial A}{\partial \vartheta}
\right).
\]
Since \(A\) is Hurwitz, this is again a Lyapunov equation, and its solution admits
the integral representation
\[
\frac{\partial X}{\partial \vartheta}
=
\int_{0}^{\infty}
e^{A^{T} t}
\left(
\frac{\partial A^{T}}{\partial \vartheta} X
+ X \frac{\partial A}{\partial \vartheta}
\right)
e^{A t}\,\mathrm{d}t .
\]
Substituting this expression into the derivative of \(f\), we obtain
\[
\frac{\partial f}{\partial \vartheta}
=
\int_{0}^{\infty}
\operatorname{tr}\!\left(
Z\,e^{A^{T} t}
\left(
\frac{\partial A^{T}}{\partial \vartheta} X
+ X \frac{\partial A}{\partial \vartheta}
\right)
e^{A t}
\right)\mathrm{d}t .
\]
By the cyclic invariance of the trace, the derivative can be written as
\[
\frac{\partial f}{\partial \vartheta}
=
\int_{0}^{\infty}
\operatorname{tr}\!\left(
e^{A t} Z e^{A^{T} t}
\left(
X\,\frac{\partial A}{\partial \vartheta}
+ \frac{\partial A}{\partial \vartheta}\,X
\right)
\right)\mathrm{d}t .
\]

We now introduce the matrix
\[
Y := \int_0^\infty e^{A t}\,Z\,e^{A^{T} t}\,\mathrm{d}t ,
\]
which is the unique solution of the dual Lyapunov equation
\[
A Y + Y A^{T} = -Z ,
\]
see, for instance, \cite{Golub1996,JBILOU2006344}.

Substituting this representation into the previous expression yields
\begin{equation}
\frac{\partial f}{\partial \vartheta}
=
\operatorname{tr}\!\left(
Y\left(
X\,\frac{\partial A}{\partial \vartheta}
+ \frac{\partial A}{\partial \vartheta}\,X
\right)
\right).
\label{eq:grad_general}
\end{equation}

\medskip

\paragraph{Derivative with respect to \(v_i\).}
For the viscosity parameter \(v_i\),
\[
\frac{\partial A}{\partial v_i} = -d(p_i)d(p_i)^{T} = -d_i d_i^{T}\,,
\]
where
\[
d_i:= d(p_i) := \begin{bmatrix} 0 \\ \Phi^T \hat{d}{(p_i)} \end{bmatrix},
\]
is a vector of dimension $2 n$, with \({\hat{d}}\)
defined in~\eqref{eq:def_hatd}.

Substituting this expression into \eqref{eq:grad_general} yields
\[
\frac{\partial f}{\partial v_i}
= -\operatorname{tr}\!\left(
Y\bigl(X d_i d_i^{T} + d_i d_i^{T} X\bigr)
\right),
\]
and by cyclicity of the trace, \(\operatorname{tr}(ABC)=\operatorname{tr}(BCA)\),
we obtain
\begin{equation}
\boxed{
\frac{\partial f}{\partial v_i}
= - \operatorname{tr}(d_i^{T} Y X d_i)
  - \operatorname{tr}(d_i^{T} X Y d_i).
} \label{derfonvi}
\end{equation}

\medskip

\paragraph{Derivative with respect to \(p_i\).}
For the position  parameter \(p_i\), the system matrix \(A\) depends on \(p_i\)
only through the vector \(d_i = d(p_i)\).
Differentiating \(A\) with respect to \(p_i\) therefore yields
\[
\frac{\partial A}{\partial p_i}
=
- \left(
w_i d_i^{T} + d_i w_i^{T}
\right),
\qquad
w_i := v_i {d'(p_i)}.
\]

The general gradient formula \eqref{eq:grad_general} gives
\[
\frac{\partial f}{\partial p_i}
=
- \operatorname{tr}\!\left(
Y\left( X (w_i d_i^{T} + d_i w_i^{T}) + (w_i d_i^{T} + d_i w_i^{T}) X \right)
\right).
\]

Expanding the trace terms and using cyclic invariance of the trace, we obtain
\begin{equation}
\boxed{
\frac{\partial f}{\partial p_i}
=
-
\Bigl(
d_i^{T} Y X w_i
+ w_i^{T} Y X d_i
+ d_i^{T} X Y w_i
+ w_i^{T} X Y d_i
\Bigr).
}
\label{derfonpi}
\end{equation}

Here \(d_i = d(p_i)\) and {$d'(p_i)$} 
are explicitly given by the finite element shape functions and their
spatial derivatives, {see Remark~\ref{rsmooth}}.

\medskip


The above derivation shows that, after computing the primal and dual
Lyapunov solutions \(X\) and \(Y\), all gradient components can be evaluated
by low-rank trace expressions, without solving additional Lyapunov equations.

\paragraph{Unified formulation of optimization criteria.}

All three optimization criteria admit a unified gradient formulation.
In each case, the objective function can be written in the trace form
\[
f(\mathbf p,\mathbf v)=\operatorname{tr}\!\bigl(Z\,X(\mathbf p,\mathbf v)\bigr),
\]
where \(X(\mathbf p,\mathbf v)\) is the unique solution of the primal Lyapunov equation
\eqref{eq:Lyap_grad_main0},
\[
A(\mathbf p,\mathbf v)^{T}X+X A(\mathbf p,\mathbf v)=-R_{hs} \,,
\]
where
\[
R_{hs} = I
\]
for the Average Total Energy criterion
(Section~\ref{subsec_average_total_energy}) and the Energy-based optimization
criterion (Section~\ref{subsec_Energy-based optimization}), and
\[
R_{hs} = Z_{ds}
\]
for the Average Total Displacements criterion (Section~\ref{subsec_average_total_displacements}).

The dual Lyapunov equation depends on the chosen optimization criterion
through the selection of the weighting matrix \(Z\), which appears as the
right--hand side of the dual equation
\[
A(\mathbf p,\mathbf v)Y+Y A(\mathbf p,\mathbf v)^{T}=-Z .
\]

More precisely, for the average total energy criterion, section \ref{subsec_average_total_energy}
 corresponding to the minimization of the mean energy over all admissible initial conditions, the weighting matrix is
\begin{equation}
\label{Def:Z0}
Z = Z_{s_0}\oplus Z_{s_0},
\qquad
Z_{s_0}=\operatorname{diag}(I_{s_0},0),
\end{equation}
where \(Z_{s_0}\) is the orthogonal projector onto the subspace spanned by the first
\(s_0\) eigenmodes; see~\cite{Nakic02, Brabender98}.

For the average total displacement criterion, section \ref{subsec_average_total_displacements}
the weighting matrix is
\[
Z_{sh} = I \,.
\]

Finally, for the energy associated with a fixed initial condition \(y_0\), section \ref{subsec_Energy-based optimization}
the objective
function is given by \(g(\mathbf p,\mathbf v)=y_0^{*}X(\mathbf p,\mathbf v)y_0\), which
corresponds to the rank--one weighting matrix
\[
Z_E = y_0 y_0^{*}.
\]
Consequently, once the primal Lyapunov solution \(X\) has been computed, the evaluation of
the gradient for all three optimization criteria differs only in the computation of the
dual Lyapunov solution \(Y\), while the subsequent gradient formulas remain identical.

\begin{algorithm}[H]
\caption{Gradient method for minimizing \( f(\pbb,\mathbf{v})=\operatorname{tr}(ZX(\pbb,\mathbf v)) \)}
\begin{algorithmic}[1]
\Require
System dimension \(n\);
number of dampers \(r\);
criterion selector \(\mathcal C\in\{\text{ATE},\text{ATD},\text{E}(y_0)\}\);
integer \(s_0\) (only if \(\mathcal C=\text{ATE}\));
initial parameters \(\pbb^{(0)} \in (0,1)^r\), \(\mathbf{v}^{(0)} \in \mathbb{R}_+^r\);
tolerance \(\varepsilon\);
maximum number of iterations \(k_{\max}\);
initial step size \(\alpha_0\);
backtracking parameters \(\beta,\sigma \in (0,1)\).

\Statex
\State \textbf{Set right-hand sides \((R_{hs},Z)\) depending on \(\mathcal C\):}
\If{\(\mathcal C=\text{ATE}\)}
    \State \(R_{hs}\gets I\)
    \State \(Z\gets Z_{s_0}\oplus Z_{s_0}\), \quad \(Z_{s_0}=\operatorname{diag}(I_{s_0},0)\)
\ElsIf{\(\mathcal C=\text{ATD}\)}
    \State \(R_{hs}\gets Z_{ds}=\begin{bmatrix}\widehat K^{-1}&0\\0&0\end{bmatrix}\)
    \State \(Z\gets I\)
\ElsIf{\(\mathcal C=\text{E}(y_0)\)}
    \State \(R_{hs}\gets I\)
    \State \(Z\gets y_0y_0^{*}\)
\EndIf

\Statex
\State \( \pbb \gets \pbb^{(0)} \), \(\mathbf{v} \gets \mathbf{v}^{(0)}\)

\For{ \(k = 0,1,\dots,k_{\max}\) } 

    \State Assemble \( A(\pbb,\mathbf{v}) \)

    \State Solve the primal Lyapunov equation for \(X^{(k)}\):
    \[
        A^T X + X A = -R_{hs}
    \]

    \State Solve the dual Lyapunov equation for \(Y^{(k)}\):
    \[
        A Y + Y A^T = -Z
    \]

    \State Compute objective value:
    \[
        f^{(k)} = \operatorname{tr}(Z\,X^{(k)}).
    \]
    \Statex
    \State \textbf{Compute gradient components}

    \For{ \( i = 1,\dots,r \) }

        \State Set \(d_i=d(p_i)\), \quad \(w_i=v_i\,d'(p_i)\).

        \State Compute
        \[
        \frac{\partial f}{\partial v_i}
        = - \operatorname{tr}\!\left( d_i^{T} Y^{(k)} X^{(k)} d_i \right)
          - \operatorname{tr}\!\left( d_i^{T} X^{(k)} Y^{(k)} d_i \right).
        \]

        \State Compute
        \[
        \frac{\partial f}{\partial p_i}
        = -\left(
            d_i^T Y^{(k)} X^{(k)} w_i
          + w_i^T Y^{(k)} X^{(k)} d_i
          + d_i^T X^{(k)} Y^{(k)} w_i
          + w_i^T X^{(k)} Y^{(k)} d_i
        \right).
        \]

    \EndFor

\algstore{GradAlg}
\end{algorithmic}
\end{algorithm}

\begin{algorithm}[H]
\ContinuedFloat
\begin{algorithmic}[1]
\algrestore{GradAlg}

 \State Form the full gradient vector
    \[
        \nabla f^{(k)} =
        \begin{bmatrix}
            \nabla_{\pbb} f^{(k)} \\[1mm]
            \nabla_{\mathbf{v}} f^{(k)}
        \end{bmatrix}.
    \]

    \If{ \( \| \nabla f^{(k)} \|_2 \le \varepsilon \) }
        \State \Return \( (\pbb,\mathbf{v}) \)
    \EndIf

    \Statex
    \State \textbf{Backtracking line search}:
    \State \( \alpha \gets \alpha_0 \)

    \While{
    \(\pbb-\alpha\nabla_{\pbb}f^{(k)}\notin(0,1)^r\) \textbf{or}
    \(\mathbf v-\alpha\nabla_{\mathbf v}f^{(k)}\not>0\) \textbf{or}
    \[
    f(\pbb - \alpha \nabla_{\pbb} f^{(k)},
      \mathbf{v} - \alpha \nabla_{\mathbf{v}} f^{(k)})
    >
    f^{(k)} - \sigma \alpha \| \nabla f^{(k)} \|^2
    \]
    }
        \State \( \alpha \gets \beta \alpha \)
    \EndWhile

    \State \textbf{Parameter update}:
    \[
        \pbb \gets \pbb - \alpha \nabla_{\pbb} f^{(k)},\qquad
        \mathbf{v} \gets \mathbf{v} - \alpha \nabla_{\mathbf{v}} f^{(k)}.
    \]

\EndFor

\end{algorithmic}
\end{algorithm}

\section{Initial Guess}

As will be shown in Section~\emph{Numerical Examples}, the optimization of damper positions on the string leads to a highly non-convex objective function with several local minima. This behaviour is fully consistent with the recent analysis of Truhar and Veselić~\cite{TruharVeselic2025} and also aligns with observations from the dimension–reduction damping optimization framework of Benner, Tomljanović and Truhar~\cite{BennerTomljTruh11}, where similarly non-convex behaviour of the objective function was reported.

Earlier studies have already indicated the inherent multi-modality of damping placement problems. In particular, Kanno et al.~\cite{KannoMatea19} report that the objective function for simultaneous optimization of positions and damping intensities is non-convex and sensitive to initialization, while G{\"u}rg{\"o}ze and M{\"u}ller~\cite{GurgozeMueller92} demonstrate that modal interactions in multi-body systems naturally generate several
distinct locally optimal damper configurations. Similar observations appear in the broader literature on supplemental damping devices, where Takewaki~(1999) and Rana--Soong~(1998) show that even simple placement or tuning tasks lead to performance surfaces with multiple stationary points. Our FEM-based formulation reproduces these phenomena, confirming that global or multi-start strategies are required to obtain reliable solutions.

\subsection{Heuristic Initial Guess}

To reduce the sensitivity of the optimization procedure to poor initialization, we develop
a heuristic strategy based on a modified explicit formula for the trace of the solution of the Lyapunov equation
\begin{equation}\label{LyapEq1Dnum-heur}
A^T X + XA = -I \,,
\end{equation}
which was originally proposed by Veseli\'c in \cite{VES90} and later presented in a refined form in \cite{Brabender98}.

For a single damper, the system matrix has the structure
\begin{equation}\label{DefMatrixAheur}
A =
\begin{bmatrix}
0 & \Omega \\
-\Omega & - v\, d d^T
\end{bmatrix},
\end{equation}
where
\[
d = (d_1,\ldots,d_n)^T \in \mathbb{R}^n
\]
is the \emph{modal damping vector} corresponding to the chosen damper position (as obtained
in Section~3).
To avoid conflict with the notation used for multiple dampers elsewhere in the paper,
we denote by
\[
\delta_k := (d)_k , \qquad k = 1,\ldots,n,
\]
the individual components of the modal damping vector associated with the single-damper
Veselić model.

For the matrix $A$ defined in \eqref{DefMatrixAheur}, and for any diagonal weighting matrix
\[
Z_\Delta = \operatorname{diag}(z_1,\ldots,z_n,z_1,\ldots,z_n),
\]
and for $X$  a solution of \eqref{LyapEq1Dnum-heur}, Veseli\'c’s identity  yields
\begin{equation}\label{TraceXZ-heur}
\operatorname{tr}(Z_\Delta X) = \frac{a}{v} + b v ,
\end{equation}
where
\begin{align}
a &= \sum_{k=s_1}^{s_0} \frac{2 z_k}{\delta_k^2},
\label{CoeffAheur} \\[0.4em]
b &= \sum_{k=s_1}^{s_0}
\left[
\frac{z_k\, \delta_k^2}{2\omega_k^2}
+
z_k \left(
 \sum_{j\neq k}
 \frac{
 3\omega_k^2 \delta_j^2 + \omega_k^2 \delta_k^2
 + \omega_j^2 \delta_j^2 + \omega_j^2 \delta_k^2
 }{(\omega_k^2 - \omega_j^2)^2}
 +
 2\,\frac{\omega_k^2}{\delta_k^2}
 \left( \sum_{j\neq k}\frac{\delta_j^2}{\omega_k^2 - \omega_j^2} \right)^2
\right)
\right]\,,
\label{CoeffBheur}
\end{align}
where $s_1$ and $s_0$ define the index range corresponding to the undamped eigenfrequencies of interest.
More precisely, for $s_1 > 1$, the eigenfrequencies of interest are
$\omega_{s_1}, \omega_{s_1+1}, \dots, \omega_{s_0}$.
This implies that the matrix $Z_\Delta$ associated with the average total energy is given by
\[
Z_\Delta = Z_{s_{1,0}} \oplus Z_{s_{1,0}},
\qquad
Z_{s_{1,0}} = \operatorname{diag}\bigl( 0_{s_1}, I_{s_0 - s_1}, 0 \bigr).
\]
If $s_1 = 1$, then $Z_\Delta$ coincides with the matrix defined in \eqref{Def:Z0}.
For more details see~\cite{Nakic02, Brabender98}.

The optimal viscosity follows explicitly:
\begin{equation}\label{vopt-heur}
v_{\mathrm{opt}} = \sqrt{\frac{a}{b}},
\end{equation}
and the corresponding optimal trace is
\begin{equation}\label{OptimalTrace-heur}
\operatorname{tr}(X_{\mathrm{opt}}) = 2 \sqrt{ab} \,,
\end{equation}
where $X_{\mathrm{opt}}$ denote the optimal solution corresponding to $Z_\Delta$.

Since the coefficients $a$ and $b$ in \eqref{CoeffAheur}–\eqref{CoeffBheur} can be
computed at a cost of $\mathcal{O}(s\, n)$, where $s=s_{0}-s_1$ denotes the number of
dominant modes appearing in the weighting matrix $Z_{\Delta}$ (typically very small),
the evaluation of the optimal trace \eqref{OptimalTrace-heur} is extremely cheap.
This allows us to sweep through all $n$ admissible damper positions and compute the
corresponding optimal values $\operatorname{tr}(X_{\mathrm{opt}}(s))$. The resulting discrete function
of the position exhibits several pronounced local minima, each of which represents a
promising candidate for the subsequent nonlinear optimization. These locally optimal
positions (together with their analytically obtained viscosities $v_{\mathrm{opt}}(s)$)
are then used as high–quality initial guesses for the gradient–based methods developed
in Sections~5 and~6.

In the present heuristic setting, the number $s$ of dominant frequencies
retained in $Z_{\Delta}$ is typically chosen between $1$ and $10$. This choice
is motivated by the empirical observation that, depending on the structural
properties of the string, the sweep usually produces between roughly $2$ and
$30$ candidate positions. These correspond to promising potential local minima
of the underlying multivariable objective and therefore provide suitable
initial guesses for the subsequent nonlinear optimization.

The resulting discrete function of the position thus exhibits several pronounced
local minima, each representing a promising candidate for the gradient--based
methods developed in Sections~5 and~6. Even for system dimensions $n < 200$,
a number of candidates exceeding $50$ would already be unnecessarily large;
in such a case, a simple uniform subdivision of the spatial domain would offer
a comparable level of coverage. Hence, restricting $s_{0}$ to a moderate range
achieves an effective balance between computational efficiency and a sufficiently
rich set of high--quality starting points.

To illustrate the observed properties, we consider the case of a single pointwise damper
applied to the damped wave equation \eqref{strong} on the spatial domain $[0,1]$,
with clamped boundary conditions \eqref{bc}, as introduced in Section~\ref{sec:introduction_preliminaries}.
The coefficients $\rho$, $A$, and $T$ are assumed to satisfy the uniform positivity conditions \eqref{positivity}, and
\[
L = 1, \qquad \rho(x) = 1, \qquad A(x) \equiv 1, \qquad T(x) \equiv 1\,.
\]

In this example, we restrict attention to a finite element discretization with $N=300$.
We present two representative cases: one with $s_{0}=10$ (where $s_1=1$ in \eqref{CoeffAheur} and \eqref{CoeffBheur}),
corresponding to a moderate number of dominant modes, and another with $s_{0}=50$ ($s_1=1$ in \eqref{CoeffAheur} and \eqref{CoeffBheur}),
illustrating the effect of a
significantly richer spectral selection. For each case, we display the corresponding
optimal damping profiles obtained from the heuristic sweep.
A more detailed analysis of this example for the case of a non-homogeneous string
(i.e., $\rho \not\equiv \mathrm{const}$) will be presented in the following section.

\begin{figure}[h!]
\centering
\includegraphics[width=0.7\textwidth]{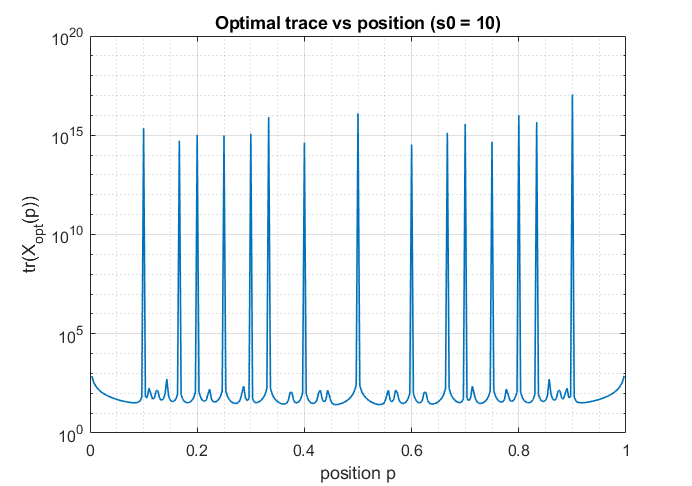}
\caption{Semilogarithmic plot of $\operatorname{tr}(X_{\mathrm{opt}}(p))$ for $s_{0}=10$.}
\label{fig:opttrace_s10}
\end{figure}

Figure~\ref{fig:opttrace_s10} displays the semilogarithmic dependence of
$\operatorname{tr}(X_{\mathrm{opt}}(p))$ on the damper position $p$, obtained
from the heuristic sweep with $s_{0}=10$. The horizontal axis represents
admissible positions $p \in (0,1)$, while the vertical axis shows the corresponding
optimal averaged energy on a logarithmic scale. Several pronounced local minima
are clearly visible, each representing a candidate position for subsequent
gradient--based optimization. From the figure, it is evident that in this case the
function exhibits $32$ local minima.

\begin{figure}[h!]
\centering
\includegraphics[width=0.7\textwidth]{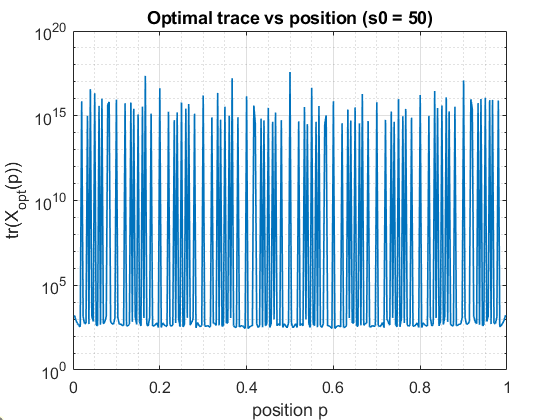}
\caption{Semilogarithmic plot of $\operatorname{tr}(X_{\mathrm{opt}}(p))$ for $s_{0}=50$.}
\label{fig:opttrace_s50}
\end{figure}

From the figure \ref{fig:opttrace_s50} corresponding to $s_0=50$, it is difficult to determine the exact number of local minima; however,
a careful inspection reveals that in this case the function exhibits approximately $100$ local minima.

If a damper is placed at a node of one of the first $s_0$ eigenfunctions, it does not dissipate the energy of
that particular eigenfunction at all. Consequently, the nodes of the first $s_0$ eigenfunctions correspond
to points of infinite energy. The local minima are therefore located in the intervals between these nodes,
which also provides a way to estimate a lower bound on the number of local minima.

\section*{Numerical examples}

As a first example in the section on numerical experiments, illustrating the properties
discussed in the previous section, we consider the case of a single pointwise damper.
In the absence of internal damping, the optimal trace is given by the
explicit formula \eqref{OptimalTrace-heur}, whose evaluation costs only
$\mathcal{O}(n\,s)$ operations, where $n$ is the model dimension and
$s=s_0-s_1$ the number of dominant frequencies.

This substantially simplifies the numerical analysis and allows us to
directly study how the parameters $s$ and $n$ influence the optimization
results. In particular, within the string model considered here, it enables
us to examine a frequently asked and still not completely resolved question:
how many dominant frequencies $s$ are needed, and which model dimension $n$
is sufficient to obtain reliable results.

In this framework, we therefore consider the damped wave equation \eqref{strong} on the spatial domain $[0,L]$,
together with the clamped boundary conditions \eqref{bc},
as introduced in Section~\ref{sec:introduction_preliminaries}.
The coefficients $\rho$, $A$, and $T$ are assumed to satisfy the
uniform positivity conditions \eqref{positivity}, i.e.,
$\rho A,\, T \in L^\infty(0,L)$ and there exist constants
$\alpha_{\rho A}>0$ and $\alpha_T>0$ such that
\[
\inf \rho A([0,L]) > \alpha_{\rho A},
\qquad
\inf T([0,L]) > \alpha_T .
\]

In this example we specify the parameters as follows:
\[
L = 1,
\qquad
T(x) \equiv 1,
\qquad
A(x) \equiv 1,
\]
so that the inhomogeneity of the string is introduced exclusively
through the density function $\rho(x)$.


We consider two density profiles.

\paragraph{(i) Borg string (\cite{Amore2010Borg}).}
Following \cite{Amore2010Borg}, the density is given by
\[
\rho(x)
=
\frac{(1+\alpha_\rho)^2}{(1+\alpha_\rho x)^4},
\qquad
\alpha_\rho = 2.
\]

\paragraph{(ii) Polynomial density string (\cite{Amore2010Borg}).}
As a second example we take
\[
\rho(x) = (1+\alpha_\rho x)^2,
\qquad
\alpha_\rho = 2.
\]

\subsection*{Damping Model}

We consider the string model and FEM setting from
Section~\ref{sec:introduction_preliminaries}, and specialize it here
to the single--damper case $r=1$ without internal damping, i.e.,
$c_0 \equiv 0$ and
\[
c = v\,\delta_p,
\qquad
p \in (0,L), \ \ v>0.
\]

\subsection*{Finite Element Discretization}

We use the standard linear finite element discretization described in
Section~\ref{sec:introduction_preliminaries} on a uniform mesh with
step size $h=L/n$, which yields the reduced matrices
\[
M_{\mathrm{red}}, \, K_{\mathrm{red}}
\in \mathbb{R}^{(n-1)\times(n-1)}.
\]
We solve the generalized eigenvalue problem
\[
K_{\mathrm{red}} \Phi
=
M_{\mathrm{red}} \Phi \Lambda,
\qquad
\Omega = \sqrt{\Lambda},
\]
and represent the localized damper in modal coordinates via the vector
$d(p)$ defined in \eqref{eq:def_hatd}--\eqref{eq:D-total}.

\subsection*{Optimization Criterion}

In this example, we consider the first two optimization criteria
introduced in Sections~\ref{subsec_average_total_energy}--\ref{subsec_average_total_displacements}.
We minimize the trace objective within the unified Lyapunov framework,
\[
f(p,v) = \operatorname{tr}\!\bigl(Z\,X(p,v)\bigr),
\]
where $X(p,v)$ is the solution of the corresponding primal Lyapunov
equation with system matrix $A(p,v)$ defined in \eqref{eq:PhaseSpaceSystem}.

As a first illustration, we sweep through all discrete damper positions.
For each position, the optimal viscosity and the corresponding optimal
trace are computed using the explicit formulas
\eqref{CoeffAheur}--\eqref{OptimalTrace-heur},
based on the first $s$ dominant eigenfrequencies.

The parameter $s$ varies from $1$ to $n-10$. For each fixed $s$, we select
the position that minimizes the trace and record the associated optimal
viscosity.

Our experiments indicate that the optimal position and the corresponding
optimal viscosity stabilize rapidly as $s$ increases. In particular,
already for $s \approx 20$, and up to $s = \mathcal{O}(0.6-0.8) \cdot n $, both quantities remain
essentially unchanged within numerical accuracy.

For the {Borg string}, the optimal position converges
to approximately $p_{\mathrm{opt}} \approx 0.2922$ with the corresponding
viscosity $v_{\mathrm{opt}} \approx 4.1537$.
For the {polynomial density} string , the optimal
position tends to $p_{\mathrm{opt}} \approx 0.6$ and the viscosity to
$v_{\mathrm{opt}} \approx 7.4$.

Therefore, in what follows, we let the matrix $Z$ select the first $s$ dominant modes;
in all experiments we set $s = 40$.

To illustrate and compare the performance of different computational
strategies, we compute the optimal position and the corresponding
optimal viscosity in three settings: the full discrete search with
$n=2000$,  the frequency cut-off model as in \cite{VESBRABDEL01,Brabender98}
with $n_{\mathrm{cut}}=200$, and  
a gradient-based method (implemented, for example, through MATLAB's
\texttt{fmincon} function).

The obtained optimal results are summarized below.
For each density profile and each computational strategy,
we report the optimal position $p^\ast$, the corresponding optimal
viscosity $v^\ast$, and the computational time.

\noindent\textbf{Borg string}

\begin{enumerate}
\item Full discrete search ($n=2000$):
\[
(p^\ast, v^\ast) = (0.2585,\; 5.9074),
\qquad
\operatorname{tr}(X^\ast) = 180,
\qquad
\text{CPU time} = 3.9\ \text{sec}.
\]

\item Frequency cut-off ($n_{\mathrm{cut}}=200$, $n=2000$):
\[
(p^\ast, v^\ast) = (0.2585,\; 5.8813),
\qquad
\operatorname{tr}(X^\ast) = 181,
\qquad
\text{CPU time} = 1.1\ \text{sec}.
\]

\item Continuous optimization using \texttt{fmincon} (with $n=100$):
\[
(p^\ast, v^\ast) = (0.292200947 ,\;   4.1537),
\qquad
\operatorname{tr}(X^\ast) = 129.4,
\qquad
\text{CPU time} = 0.66\ \text{sec}.
\]
\end{enumerate}

\medskip
\noindent\textbf{Polynomial density string}

\begin{enumerate}
\item Full discrete search ($n=2000$):
\[
(p^\ast, v^\ast) = (0.62,\; 8.7676),
\qquad
\operatorname{tr}(X^\ast) = 294.6,
\qquad
\text{CPU time} = 3.9562\ \text{sec}.
\]

\item Frequency cut-off ($n_{\mathrm{cut}}=200$, $n=2000$):
\[
(p^\ast, v^\ast) = (0.62,\; 8.7372),
\qquad
\operatorname{tr}(X^\ast) = 309.96,
\qquad
\text{CPU time} = 1.07\ \text{sec}.
\]

\item Continuous optimization using \texttt{fmincon} (with $n=100$):
\[
(p^\ast, v^\ast) = (0.599411615  ,\;   0.747962829),
\qquad
\operatorname{tr}(X^\ast) = 251.47,
\qquad
\text{CPU time} = 0.047\ \text{sec}.
\]
\end{enumerate}

\noindent
Observe that the \texttt{fmincon} approach yields a smaller value of
$\operatorname{tr}(X^\ast)$ in both examples. However, since the
optimization is performed with $n=99$, this discretization is too
coarse for $s=40$ to accurately resolve the relevant frequency range,
which explains the discrepancy with the results obtained for
$n \gtrsim 400$.  For example, for the Borg string with $n=999$ and
the starting point $(0.2585,\,5.9074)$, the optimization yields
$(p^\ast, v^\ast) = (0.25888797, \, 5.7481297)$,
with the corresponding optimal trace
$\operatorname{tr}(X^\ast) = 175.054686$.
For the Polynomial density string with $n=999$ and
the starting point $(0.625 , \, 8)$, the optimization yields
$(p^\ast, v^\ast) = (0.627533322045556, \,   8.815503737293213)$,
with the corresponding optimal trace
$\operatorname{tr}(X^\ast) =  310.40528$.

\subsection*{Illustration of the sufficient model dimension and number of dominant frequencies.}

In the following example, we illustrate the practical choice of the
model dimension $n$ (with $n = N-1$) and the number of dominant
frequencies $s$ required to reliably determine the optimal damper
position and the corresponding optimal viscosity for the average energy criterion.
We consider two pointwise dampers at $p_1, p_2$ with viscosities
$v_1, v_2$, and include a small uniform internal damping
$c_0 = 10^{-7}$.

To this end, we consider a homogeneous string, which enables a direct
comparison between the discrete approximation and the exact undamped spectrum.
In the homogeneous case, all coefficients are constant, namely
\[
\rho(x) = 1,
\qquad
T(x) = 1,
\qquad
A(x)=1,
\qquad
L = 1.
\]
The exact undamped eigenfrequencies are therefore given by
\[
\omega_k = k\pi, \qquad k \in \mathbb{N}.
\]

We consider three discretization levels,
\[
N = 200, \quad 1000, \quad 2000,
\]
which correspond to model dimensions $n = N-1$.
The relative error of the $t$-dimensional approximation of the
$s$-th eigenfrequency is defined as
\[
rerr_t(s)
=
\frac{|\omega_t(s) - s\pi|}{s \pi},
\qquad
t \in \{200,1000,2000\}.
\]

For $s = 40$, we obtain
\[
rerr_{200}(40)=0.0165 < 0.1,
\]
while the corresponding errors
$rerr_{1000}(40)$ and $rerr_{2000}(40)$
are significantly smaller.
We therefore compare the optimal damper positions,
the corresponding optimal viscosities,
and the resulting optimal trace values for $s = 40$
across all three model dimensions $N=200,1000,2000$.

The computed optimal damper positions $p_{\mathrm{opt}}$, the corresponding
optimal viscosities $v_{\mathrm{opt}}$, and the resulting optimal trace
values for $s=40$ are summarized in
Table~\ref{tab:homogeneous-convergence}. The results clearly indicate
convergence with respect to the model dimension $n$.

\begin{table}[h!]
\centering
\caption{Optimal damper positions, viscosities, and trace values
for $s=40$ at different discretization levels.}
\label{tab:homogeneous-convergence}
\begin{tabular}{c|cc|cc|c}
\hline
$N$
& $p_{\mathrm{opt}}^{(1)}$ & $p_{\mathrm{opt}}^{(2)}$
& $v_{\mathrm{opt}}^{(1)}$ & $v_{\mathrm{opt}}^{(2)}$
& $\mathrm{Trace}_{\mathrm{opt}}$ \\
\hline
200  & 0.2661 & 0.7425 & 2.3663 & 2.4265 & 57.6751 \\
1000 & 0.2208 & 0.7333 & 2.0865 & 2.0966 & 57.5660 \\
2000 & 0.2208 & 0.7332 & 2.0823 & 2.0938 & 57.6137 \\
\hline
\end{tabular}
\end{table}

The candidate positions for local minima obtained from \eqref{OptimalTrace-heur} amount
to $32$ (considering only half of the domain due to the symmetry of the string).
Among these, the candidates indexed from the $13$th to the $17$th are given by
\[
0.1960,\quad 0.2211,\quad 0.2312,\quad 0.2462,\quad 0.2563.
\]
In particular, the candidates between the $15$th and $17$th positions appear to be closest
to the true local minima.

In addition, for the finest discretization $N = 2000$,
we also compute the optimal configuration for $s = 200$.
Since in this case
\[
rerr_{2000}(200) = 0.0041 < 0.1,
\]
this comparison illustrates how far into the spectrum one needs to proceed
in order to obtain a practically sufficient and spectrally reliable solution
of the optimization problem.

For $N = 500$, we examine the dependence of the optimal configuration
on the number of dominant frequencies $s$.
The computed optimal damper positions
$p_{\mathrm{opt}}^{(1)}$, $p_{\mathrm{opt}}^{(2)}$,
the corresponding optimal viscosities
$v_{\mathrm{opt}}^{(1)}$, $v_{\mathrm{opt}}^{(2)}$,
and the resulting optimal trace values
for $s = 20,40,60,80,100$ are summarized in
Table~\ref{tab:homogeneous-s-dependence}.

The results show evident differences in the optimal damper positions
as $s$ increases; however, these variations occur predominantly in the
second decimal place.
Since in our model we employ $P1$ finite elements, the spatial
approximation is of order $\mathcal{O}(h)$ in the $H^1$–norm.
For $N=500$, the mesh size satisfies $h \approx 1/N$, and therefore
the attainable accuracy in the computed positions is consistent with
an $\mathcal{O}(h)$ discretization error.
The optimal viscosities adjust accordingly as additional spectral
components are included, while the optimal trace value increases with $s$
due to the contribution of a larger number of modes in the objective functional.

\begin{table}[h!]
\centering
\caption{Optimal damper positions, viscosities, and trace values
for $N=500$ and different numbers of dominant frequencies $s$.}
\label{tab:homogeneous-s-dependence}
\begin{tabular}{c|cc|cc|c}
\hline
$s$
& $p_{\mathrm{opt}}^{(1)}$ & $p_{\mathrm{opt}}^{(2)}$
& $v_{\mathrm{opt}}^{(1)}$ & $v_{\mathrm{opt}}^{(2)}$
& $\mathrm{Trace}_{\mathrm{opt}}$ \\
\hline
20  & 0.2769 & 0.6821 & 1.9578 & 1.8762 & 27.6092 \\
40  & 0.2615 & 0.6404 & 2.2260 & 1.9786 & 59.2590 \\
60  & 0.2777 & 0.7372 & 2.0698 & 2.1090 & 89.2564 \\
80  & 0.2595 & 0.7299 & 2.1225 & 2.1295 & 120.7392 \\
100 & 0.2577 & 0.7339 & 2.1481 & 2.2121 & 153.2824 \\
\hline
\end{tabular}
\end{table}

It is important to emphasize that the obtained results correspond to local minima,
and no guarantee can be provided that they are globally optimal.
What is evident, however, is that for smaller values of $s$, the optimization process is more reliable.
This is due to the fact that the spectral approximation is more accurate for lower frequencies,
resulting in smaller errors in the computed eigenvalues, as well as to the considerably smaller
number of local minima. Consequently, the optimization procedure benefits
from increased accuracy when fewer (i.e., lower) dominant frequencies are taken into account.

\subsection*{Numerical example: two dampers for aeolian vibrations (real-life setting)}

In this example we consider a model motivated by the classical work of Hagedorn~\cite{Hagedorn1982},
where wind--excited vibrations of overhead transmission lines are studied. In particular, we focus on aeolian vibrations
caused by vortex shedding, which typically occur in the frequency range $10$--$50$~Hz and may lead to fatigue damage of the cable.

We consider a typical $220\,\mathrm{kV}$ transmission line with span length $L=300\,\mathrm{m}$ and a standard ACSR conductor,
with mass density satisfying $\rho(x)A(x)=2\,\mathrm{kg/m}$.
The mechanical tension is high (on the order of $T \approx 3.7 \times 10^4\,\mathrm{N}$) and slightly increases towards the supports.
Accordingly, we use a spatially varying tension of the form
\[
\rho(x)A(x) = 2\,\mathrm{kg/m}, \qquad
T(x)= 37000\cdot \left(1 + 0.1 \cdot \frac{|x - L/2|}{L/2}\right),
\]
which serves as a simple approximation of the increase of tension towards the supports, consistent with the catenary model.

The cable is modeled by the wave equation \eqref{strong} on the spatial domain $[0,L]$ with clamped boundary conditions \eqref{bc},
using the same FEM discretization \eqref{eq:SecondOrderFEM} and phase--space formulation \eqref{Def:PhaseSpacevec}--\eqref{eq:PhaseSpaceSystem} as in the previous sections.

\paragraph{Frequency range and interpretation.}

Aeolian vibrations are associated with wind-induced excitation frequencies typically lying in the range from approximately $5$~Hz up to about $24$~Hz, and in some practical situations extending up to about $48$~Hz.

In terms of angular frequencies, this corresponds to
\[
31 \lesssim \omega_i \lesssim 150,
\]
and, for the extended range,
\[
31 \lesssim \omega_i \lesssim 300.
\]

The natural frequencies of the considered model are obtained from the generalized eigenvalue problem and are expressed in terms of angular frequencies $\omega_i$. The corresponding frequencies in Hz are given by
\[
f_i = \frac{\omega_i}{2\pi}.
\]

In the present discretization, the computed spectrum satisfies $\omega_1 \approx 0.22$ and $\omega_{n} \approx 234.3$, which implies that the model captures frequencies up to approximately
\[
f_{\max} \approx \frac{234.3}{2\pi} \approx 37~\text{Hz}.
\]
Therefore, the considered model fully covers the range up to approximately $24$~Hz and partially covers the extended range up to approximately $48$~Hz.

Since the tension $T(x)$ is non-constant, the eigenfrequencies do not follow a simple analytical law, and no direct correspondence between modal indices and frequency ranges can be assumed. Accordingly, the relevant part of the spectrum is identified directly from the computed eigenvalues, by selecting those modes whose angular frequencies satisfy $\omega_i \gtrsim 31$.

Moreover, a simple one-dimensional approximation, based on the heuristic formula~\eqref{OptimalTrace-heur}, suggests that for frequencies in the range $\omega \in [31,150]$ (and up to $\omega \in [31,300]$), the optimal position of a single damper lies approximately between $0.5\,\mathrm{m}$ and $2\,\mathrm{m}$ from the support.

In particular, higher frequencies tend to correspond to positions closer to the support (i.e., smaller values of $p$). This trend is consistent with the qualitative description of aeolian vibrations in~\cite{Hagedorn1982} and with standard engineering practice, where Stockbridge dampers are typically installed within a few meters of the support.

\paragraph{Dampers.}

We consider configurations with $r=2$ and $r=4$ pointwise dampers (representing Stockbridge dampers),
located at positions $p_1,\dots,p_r \in (0,L)$ with corresponding viscosities $v_1,\dots,v_r>0$.

Since the resulting optimization problem exhibits a highly nontrivial structure with multiple local minima
noting that in the general case the number of local minima increases rapidly with $s$, we restrict the analysis
to these practically relevant configurations, obtained by the heuristic formula~\eqref{OptimalTrace-heur}.

\paragraph{Optimization problem.}
The goal is to determine the optimal positions and viscosities by minimizing the trace functional
\[
f(p_1,\dots,p_r,v_1,\dots,v_r)
=
\tr\!\bigl(ZX(p_1,\dots,p_r,v_1,\dots,v_r)\bigr),
\]
for $r\in\{2,4\}$, where $X$ solves the Lyapunov equation \eqref{eq:Lyap_grad_main0}.
The weighting matrix $Z$ is defined as
\[
Z = Z_1 \oplus Z_1,
\qquad
Z_1 = 0_{31} \oplus I_{120} \oplus 0_{n-151},
\]
so that the functional selects a prescribed band of dominant modes.
This corresponds to minimizing the averaged vibration energy over the relevant spectral range.

This setting is consistent with the energy-based approach of Hagedorn~\cite{Hagedorn1982},
where the envelope of the response, rather than individual resonance peaks,
is of primary importance for practical design and damper tuning.

The obtained optimal positions and corresponding viscosities for the considered
configurations are summarized in the following table.

\begin{table}[h!]
\centering
\caption{Optimal damper configurations for the considered aeolian vibration setting.}
\begin{tabular}{c|c|c}
\hline
$r$ & Positions $p_i$ & Viscosities $v_i$ \\
\hline
2
& $\begin{aligned}
   &[1.5356,\; 298.4644]
  \end{aligned}$
& $\begin{aligned}
   &[331.4824 ,\; 331.4824]
  \end{aligned}$ \\

4
& $\begin{aligned}
   &[1.5715,\; 3.1514,\\
   &\quad 298.4285,\; 296.8486]
  \end{aligned}$
& $\begin{aligned}
   &[348.6473,\; 133.0398,\\
   &\quad 348.6474,\; 133.0398]
  \end{aligned}$ \\
\hline
\end{tabular}
\end{table}
The configurations exhibit a clear symmetry with respect to the midpoint of the span.

\newpage

\subsection*{Numerical example: two dampers and three Lyapunov criteria (localized density)}

In this example, we illustrate the effects of Average Total Energy (ATE) and Average Total Displacement (ATD), alongside optimization strategies for specific initial data.
As shown in Figure \ref{fig:opttrace_s10}, for $s = 10$ with a single damper, the cost function exhibits numerous local minima;
the quantity of these minima is directly related to the number of nodes in the eigenvectors we aim to damp.
With additional dampers, the number of local minima increases significantly. Furthermore, as the values of these local minima are
nearly identical (they differ up to 10\%),
the minimization procedure struggles to isolate the global minimum.  Thus in this example we illustrate different optimality criteria for $s = 4$,
using an initial condition defined as the sum of the first four eigenvectors (with zero velocity).

For each optimality criterion, we determine the optimal positions and viscosities of the dampers.
We then plot the displacements at four distinct time points and present the time evolution of the total energy for each case.
Our example indicate that the ATE criterion damps the energy more rapidly than the ATD.
However, it should be noted that this behavior may not hold for all initial data, as both criteria are formulated based on averages.
Furthermore, dampers optimized for specific initial data show the highest performance, as they are tailored exactly to those conditions.

We remain within the framework introduced in
Sections~\ref{sec:introduction_preliminaries}--\ref{Sec:Gradoftracefunc}.
In particular, we use the same FEM model \eqref{eq:SecondOrderFEM},
the phase--space transformation \eqref{Def:PhaseSpacevec}--\eqref{eq:PhaseSpaceSystem},
and the unified trace formulation
\(
f(p_1,p_2,v_1,v_2)=\tr(ZX(p_1,p_2,v_1,v_2))
\)
with the primal and dual Lyapunov equations
\eqref{eq:Lyap_grad_main0} and
\(
A(p_1,p_2,v_1,v_2)Y+YA(p_1,p_2,v_1,v_2)^T=-Z
\),
as described in Section~\ref{Sec:Gradoftracefunc}.

\paragraph{Coefficients and density.}
We set \(L=1\), \(T(x)\equiv 1\), and \(A(x)\equiv 1\), so that the
inhomogeneity enters only through the density \(\rho(x)\).
Instead of the Borg/polynomial profiles used earlier, we take the smooth localized
profile
\begin{equation}\label{eq:rho-bump-example}
\rho(x)=10\,x(1-x)\exp\!\bigl(-40(x-0.2)^2\bigr)+0.1,\qquad x\in[0,1].
\end{equation}

This choice produces a smooth localized peak near \(x=0.2\),
while the additive constant \(0.1\) guarantees strict positivity
on the entire interval. The corresponding density profile is shown
in Figure~\ref{fig:rho-bump-example}.

\begin{figure}[H]
\centering
\includegraphics[width=0.65\textwidth]{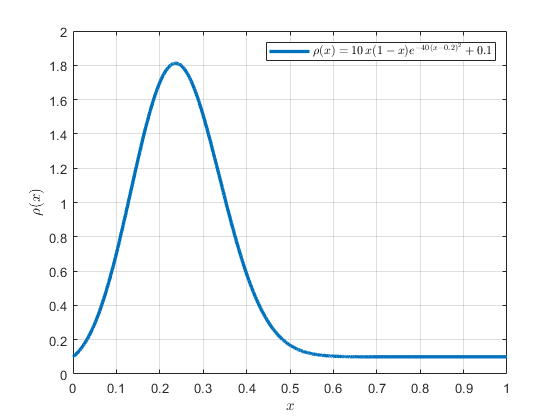} 
\caption{Smooth localized density profile \(\rho(x)\) defined in \eqref{eq:rho-bump-example}.}
\label{fig:rho-bump-example}
\end{figure}

\paragraph{Two dampers and initial data.}
We consider \(r=2\) pointwise viscous dampers without internal damping,
\(
c_0\equiv 0
\),
so that
\[
c(x)=v_1\delta_{p_1}(x)+v_2\delta_{p_2}(x),
\qquad
0<p_1<p_2<1,\quad v_1,v_2>0.
\]

\paragraph{Discretization and dominant frequencies.}
All computations are performed with \(N=200\) (hence \(n=N-1=199\)) and
we select the first \(s=4\) dominant frequencies in the average--energy criterion,
as in \eqref{Z_energy}.

\paragraph{Three optimization criteria.}
We compute optimal damper positions and viscosities for the following three
objectives, all written in the unified trace form
\(f(\pbb,\vbb)=\tr(ZX(\pbb,\vbb))\):
\begin{itemize}
\item \emph{Average total energy} (ATE), Section~\ref{subsec_average_total_energy},
      with \(R_{hs}=I\) and \(Z=Z_{s}\oplus Z_{s}\).
\item \emph{Average total displacements} (ATD), Section~\ref{subsec_average_total_displacements},
      with \(R_{hs}=Z_{ds}\) and \(Z = I\). 
\item \emph{Energy for a fixed initial state} \(E(y_0)\),
      Section~\ref{subsec_Energy-based optimization}, with \(R_{hs}=I\) and \(Z=y_0y_0^{T}\).
\end{itemize}

\paragraph{Multi-start initialization.}
To mitigate sensitivity to initialization and to detect potential local minima,
we use the multi-start strategy described in Section~\emph{Initial Guess}.
Specifically, for each criterion we initialize the nonlinear solver from the discrete set
\[
(p_1^{(0)},p_2^{(0)},v_1^{(0)},v_2^{(0)})
=
\left(\frac{k}{25},\;0.25+\frac{k}{35},\;1,\;1\right),
\qquad k=1,\dots,20,
\]
and we report the best configuration among the obtained local minimizers.

\paragraph{Results.}

The multi-start optimization revealed several (approximately ten)
distinct local minima. Among the obtained stationary configurations,
we select and report the one with the smallest objective value.

For illustration,  local minimizer for  ATE is obtained
after restarting with $x_0=[0.2, 0.54,  2, 1]$:
\[
(p_1,p_2,v_1,v_2)
=
( 0.153426,\;  0.6272614,\;  2.180858,\; 0.643924 )
\]
for which
\[
\operatorname{tr}(X_0)=2.5048.
\]
The corresponding gradient is of order $10^{-8}$,
\[
\nabla \operatorname{tr}(X_0)
=
10^{-8}\,(1.08,\, 5.23,\, 1.17,\, 5.36)^T \,,
\]
confirming that this configuration satisfies the first-order optimality
conditions up to numerical tolerance, although it does not correspond
to the minimal objective value among all detected stationary points.

The optimal configurations \((p_1^\ast,p_2^\ast,v_1^\ast,v_2^\ast)\) and the
corresponding optimal objective values are summarized in
Table~\ref{tab:two-damper-bump-three-criteria}.

\begin{table}[H]
\centering
\caption{Two-damper optimization results for the localized density \eqref{eq:rho-bump-example}
(\(N=200\), \(n=199\), \(s=4\)), $c_0=0.0001$.}
\label{tab:two-damper-bump-three-criteria}
\begin{tabular}{l|cc|cc|c}
\hline
Criterion
& $p_1^\ast$ & $p_2^\ast$
& $v_1^\ast$ & $v_2^\ast$
& Optimal value \\
\hline
Average total energy (ATE)
& 0.1534 &  0.62726
& 2.1809  &  0.6439
& 2.5048 \\
Average total displacement (ATD)
& 0.2104 & 0.6392
& 2.3303 & 0.7558
& 12.0074 \\
Fixed-initial-state energy $E(y_0)$
& 0.1210 & 0.3871
& 1.6884 & 1.9882
& 0.19105 \\
\hline
\end{tabular}
\end{table}

We illustrate the displacement profiles (solutions) for the optimal system
\(
A_{\text{opt}} = A(p_1^\ast, p_2^\ast, v_1^\ast, v_2^\ast)
\)
under the average total energy (ATE) criterion in Figure \ref{fig:two-damper-bump-ATE}.
The displacement component \( x_{p_1}(x) \) is shown at selected time instants
\(t = 0.001,\ 0.701,\ 1.401,\ 2.101\).
This provides a visual representation of how the optimally placed dampers
affect the system’s displacement over time.

\begin{figure}[H]
\centering
\includegraphics[width=0.9\textwidth]{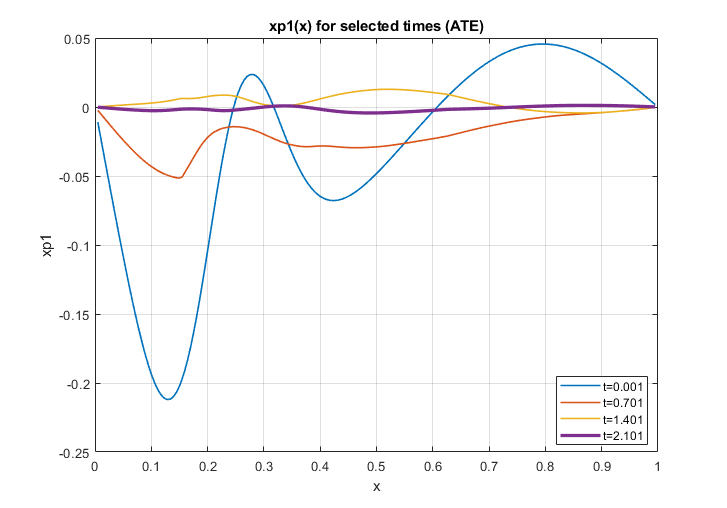}
\caption{
Spatial profile of the projected displacement component
\(
x_{p_1}(x)
=
d_0^{T}\bigl(\Omega^{-1} y(t)\bigr)
\)
corresponding to the optimal configuration
\(
(p_1^\ast,p_2^\ast,v_1^\ast,v_2^\ast)
=
(0.1534,\,0.6272,\,2.1808,\,0.6439)
\)
obtained for the \emph{average total energy} (ATE) criterion.
The curves are shown for four representative time instants
\(t = 0.001,\ 0.701,\ 1.401,\ 2.101\),
as indicated in the legend.
The thick curve corresponds to the largest time and illustrates
the long–time decay behavior induced by the optimally placed dampers.
}
\label{fig:two-damper-bump-ATE}
\end{figure}

We  present the displacement profiles (solutions) for the optimal system \\
\( A_{\text{opt}} = A(p_1^\ast, p_2^\ast, v_1^\ast, v_2^\ast) \) under the average total energy (ATD) criterion
in Figure \ref{fig:two-damper-bump-ATD}.
The displacement component \( x_{p_1}(x) \) is shown at selected time instants \\
\( t = 0.001,\ 0.701,\ 1.401,\ 2.101 \).
This illustrates how the optimally placed dampers for ATD optimizaiton influence the evolution and decay of
the system’s displacement over time.

\begin{figure}[H]
\centering
\includegraphics[width=0.8\textwidth]{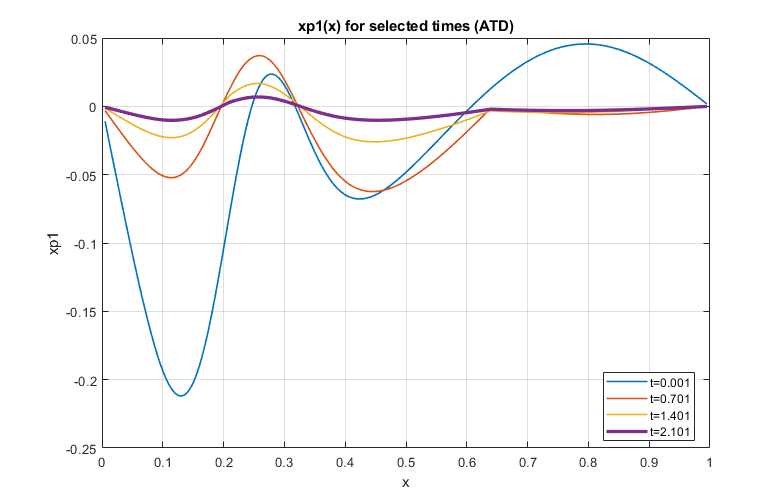}
\caption{
Spatial profile of the projected displacement component
\(
x_{p_1}(x)
=
d_0^{T}\bigl(\Omega^{-1} y(t)\bigr)
\)
corresponding to the optimal configuration
\(
(p_1^\ast,p_2^\ast,v_1^\ast,v_2^\ast)
=
(0.2104,\,0.6392,\,2.3303,\,0.7558)
\)
obtained for the \emph{average total displacement} (ATD) criterion.
The curves are shown for four representative time instants
\(t = 0.001,\ 0.701,\ 1.401,\ 2.101\),
as indicated in the legend.
The thick curve corresponds to the largest time and reflects
the spatial redistribution of displacement induced by the
optimal damper placement under the displacement-based criterion.
}
\label{fig:two-damper-bump-ATD}
\end{figure}

As the final case, we present the displacement components of the solution obtained for the optimal system
\(
A_{\text{opt}} = A(p_1^\ast, p_2^\ast, v_1^\ast, v_2^\ast)
\)
under the energy for a fixed initial state \(E(y_0)\) criterion in Figure \ref{fig:two-damper-bump-Ey0}.
The displacement component
\(
x_{p_1}(x)
\)
is shown at the time instants
\(t = 0.001,\ 0.701,\ 1.401,\ 2.101\).

\begin{figure}[H]
\centering
\includegraphics[width=0.9\textwidth]{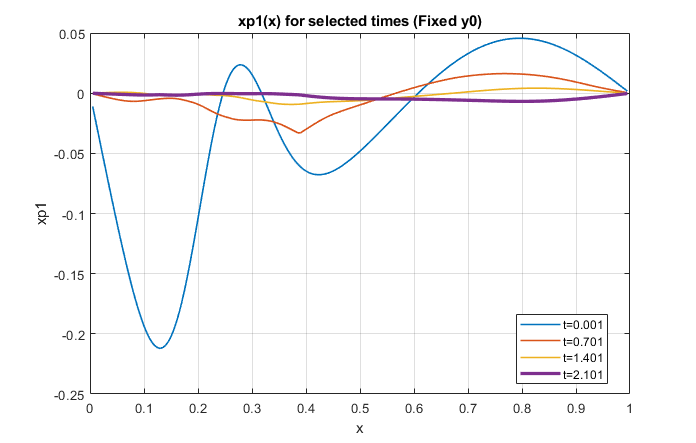}
\caption{
Spatial profile of the projected displacement component
\(
x_{p_1}(x)
=
d_0^{T}\bigl(\Omega^{-1} y(t)\bigr)
\)
corresponding to the optimal configuration
\(
(p_1^\ast,p_2^\ast,v_1^\ast,v_2^\ast)
=
(0.1300,\,0.3875,\,1.7243,\,2.0328)
\)
obtained for the \emph{energy for a fixed initial state} \(E(y_0)\) criterion.
The curves are shown for four representative time instants
\(t = 0.001,\ 0.701,\ 1.401,\ 2.101\),
as indicated in the legend.
The thick curve corresponds to the largest time and illustrates
the long--time decay behavior induced by the optimally placed dampers.
}
\label{fig:two-damper-bump-Ey0}
\end{figure}

Simultaneously, we present a comparison of the corresponding energy responses over the considered time interval
in Figure \ref{fig:energy-comparison}.

\begin{figure}[H]
\centering
\includegraphics[width=0.9\textwidth]{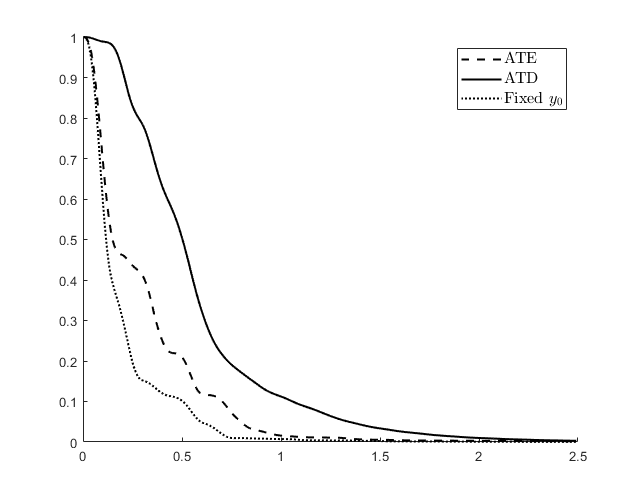}
\caption{
Comparison of the corresponding energy responses over the considered time interval for the three optimization criteria: average total energy (ATE), average total displacement (ATD), and the fixed initial state energy \(E(y_0)\).
The different line styles ensure clear distinction in grayscale representation.
}
\label{fig:energy-comparison}
\end{figure}

To further assess the influence of the initial condition, we modify the initial displacement such that
\( x_0(0) \) corresponds to the \((s+1)\)-th eigenvector, (recall $\dot{x}_0=0$).
 In this case, an additional optimal component \( x_{\mathrm{opt}}^{Y} \) is obtained, given by
\[
(p_1^\ast,p_2^\ast,v_1^\ast,v_2^\ast)_{y_0} =
\begin{bmatrix}
2.1111\cdot 10^{-1} & 7.9912\cdot 10^{-1} & 2.4563 & 3.9480\cdot 10^{-1}
\end{bmatrix}.
\]
This corresponds to the contribution associated with the modified initial state \( y_0 \),
and the resulting energy evolution reflects the increased influence of higher-order modes.

Note that in this example the optimal positions for ATE are computed for $s = 4$, i.e., only the first four modes are taken into account in determining the optimal configuration. Therefore, in contrast to the case shown in Figure \ref{fig:energy-comparison}, it is not surprising that the optimal ATD
configuration reduces the total energy faster than ATE. Moreover, the energy decay for the ATE configuration is significantly slower than in Figure~\ref{fig:energy-comparison}, for which the optimal ATE configuration is constructed.

\begin{figure}[H]
\centering
\includegraphics[width=0.9\textwidth]{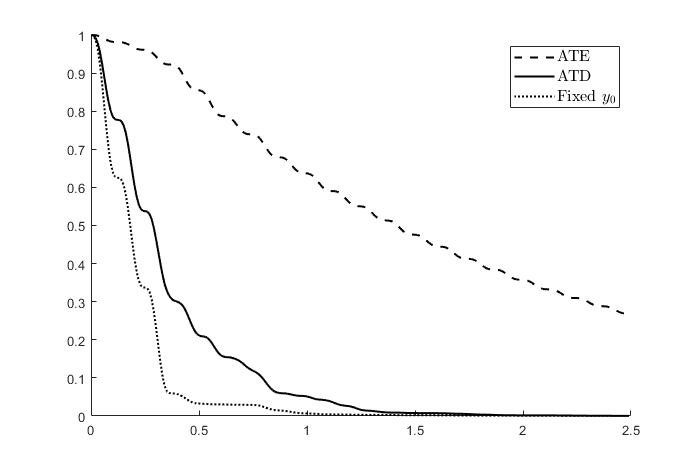}
\caption{
Comparison of the corresponding energy responses over the considered time interval for the three optimization criteria: average total energy (ATE), average total displacement (ATD), and the fixed initial state energy \(E(y_0)\).
The different line styles ensure clear distinction in grayscale representation.
}
\label{fig:energy-comparison2b}
\end{figure}

\paragraph{Comparison of ATE and ATD and Fixed Initial Condition.}

\begin{itemize}
\item The optimization problem is highly nonconvex, with many local minima of similar value, making the identification of a global minimizer challenging and sensitive to initialization.

\item The ATE criterion generally achieves faster energy decay than ATD for the considered modes; however, this advantage depends on the initial condition, and changes in the modal content can significantly alter the optimal damper configuration.

\item The fixed--initial--state criterion \(E(y_0)\) provides the best performance when the initial state is known, highlighting the importance of selecting an appropriate optimality criterion and using robust optimization strategies.
\end{itemize}

\paragraph{Conclusion}

In this paper, we have developed a framework for the optimization
of pointwise viscous dampers in vibrating string models by combining
a finite element discretization with a phase--space reformulation
amenable to Lyapunov-based analysis.

A key contribution is the unified treatment of three optimization
criteria---average total energy (ATE), average total displacement (ATD),
and energy for a fixed initial state---all expressed as trace minimization
problems involving Lyapunov equations. This leads to a single computational
pipeline in which the dominant cost is the solution of a pair of Lyapunov
equations.

We derived compact gradient formulas in terms of the primal and dual
Lyapunov solutions, yielding an efficient implementation that avoids
additional Lyapunov solves.

The resulting optimization problem is highly nonconvex, with many local
minima and strong sensitivity to initialization. To mitigate this, we
introduced a heuristic initialization based on an explicit single-damper
trace formula. This heuristic is essential for practical performance:
it delivers high-quality starting points at negligible cost and enables
the gradient method to reliably converge to meaningful local minimizers.

Numerical experiments show that (i) optimal positions and viscosities
stabilize quickly once the selected spectral band is sufficiently well
resolved, (ii) inadequate model dimensions may lead to misleading optima,
and (iii) the choice of criterion is problem-dependent: ATE favors faster
average energy decay, ATD may be preferable for certain displacement-focused
responses, while the fixed-initial-state criterion performs best when the
excitation is known.

Finally, the aeolian vibration example demonstrates the practical relevance
of the approach, with optimal configurations consistent with engineering
practice.

Overall, reliable damping optimization requires a consistent combination of
spectral accuracy, appropriate mode selection, and robust initialization.
Future work will address extensions to more complex structures, integration
with model reduction techniques, and strategies for handling the multimodal
nature of the problem.

\section*{Acknowledgments}

J. Tamba\v{c}a research has been supported by the Croatian Science Foundation under the project number HRZZ-IP-2022-10-1091, by the project "Implementation of cutting-edge research and its application as part of the Scientific Center of Excellence for Quantum and Complex Systems, and Representations of Lie Algebras", Grant No. PK.1.1.10.0004, by the European Union -- NextGenerationEU through the National Recovery and Resilience Plan 2021-2026, institutional grant of University of Zagreb Faculty of Science (IK IA 1.1.3. Impact4Math).

Ninoslav Truhar was partially supported by the Croatian Science Foundation under the project \textit{Optimal control and model reduction for evolution and data driven problems} (Conduction; Grant No. HRZZ-IP-2022-10-5191) and was partially supported by the European union-NextGenerationEU, grant no. 581-UNIOS-55, OpHoMat.

\bibliographystyle{plain}
\bibliography{bibliogr_opt_abNT}

\end{document}